\RequirePackage[english]{babel}
\documentclass[11pt]{amsart}
\usepackage{hyperref}
\usepackage{a4wide}
\usepackage{amsmath}
\usepackage{amsthm}
\newcommand{\RR}{\mathbb{R}}
\newcommand{\CC}{\mathbb{C}}
\newcommand{\NN}{\mathbb{N}}
\newcommand{\ZZ}{\mathbb{Z}}
\DeclareMathOperator{\EE}{\mathbb{E}}
\newcommand{\PP}{\mathbb{P}}
\newcommand{\cD}{\mathcal{D}}
\newcommand{\cV}{\mathcal{V}}
\newcommand{\G}{G}
\newcommand{\E}{E}

\DeclareMathOperator{\sign}{\mathrm{sign}}
\newcommand{\CG}{C^\Gamma}
\newcommand{\Tr}{{\mathop{\mathrm{Tr}}}}
\DeclareMathOperator{\vol}{vol}
\DeclareMathOperator{\supp}{\mathrm{supp}}

\newtheorem{thm}{Theorem}
\newtheorem{lem}[thm]{Lemma}
\theoremstyle{definition}
\newtheorem{dfn}[thm]{Definition}
\theoremstyle{remark}
\newtheorem{rem}[thm]{Remark}

\begin{document}
\title[Wegner estimates for Schr\"odinger operators
on metric graphs]{The modulus of continuity of Wegner estimates for random Schr\"odinger operators
on metric graphs}

\author[M.~J.~Gruber]{Michael J.\ Gruber}
\address{TU Clausthal\\
Institut f\"ur Mathematik\\
38678 Clausthal-Zellerfeld\\
Germany}

\urladdr{\url{http://www.math.tu-clausthal.de/~mjg/}}	

\author[I.~Veseli\'c]{Ivan Veseli\'c}
\address{TU Chemnitz\\
Fakult\"at f\"ur Mathematik\\
09107 Chemnitz\\
Germany}
\urladdr{\url{http://www.tu-chemnitz.de/mathematik/schroedinger/members.php}}

\thanks{The authors were financially supported by the DFG under grant   Ve 253/2-2 within the Emmy-Noether-Programme.}

\thanks{\copyright 2007 by the authors. Faithful reproduction of this article,
         in its entirety is permitted for non-commercial purposes. {\today, \jobname.tex}}


\subjclass[2000]{}

\begin{abstract}
We consider an alloy type potential on an infinite metric graph.
We assume a covering condition on the single site potentials.
For random Schr\"odingers operator associated with the alloy type potential 
restricted to  finite volume subgraphs
we prove a Wegner estimate which reproduces the modulus of continuity of 
the single site distribution measure.  The Wegner constant is independent of the energy.
\end{abstract}

\maketitle

\section{Introduction}

We study the distribution  of eigenvalues of alloy-type random Schr\"odinger 
operators on finite metric graphs. More precisely, we prove a 
Wegner estimate on the average number of eigenvalues in a given
energy interval. Our (upper) bound is proportional to the 
volume of the finite metric graph and reproduces the modulus of 
continuity of the single site random coupling constants.
The volume is the sum of the lengths of all edges. In particular, 
for equilateral graphs it is the number of edges times the single edge length.

In the case that there exists a selfaveraging integrated density of states
our result implies estimates on its modulus of continuity.

Wegner estimates have been studied for Anderson-type operators on $l^2(\ZZ^d)$
and alloy-type operators on $L^2(\RR^d)$ starting with the paper \cite{Wegner-81}.
The strategy presented here follows the line of argument of \cite{HundertmarkKNSV-06}.
There the same type of Wegner bound as in this note was proven
for discrete operators on  $l^2(\ZZ^d)$. For alloy-type operators on $L^2(\RR^d)$ 
a similar upper bound was derived. There however, a logarithmic correction term
due to the (possible) singularity of the spectral shift function appears.

In the  recent \cite{CombesHK} a Wegner estimate for
random Schr\"odinger operators on $L^2(\RR^d)$ without the logarithmic 
correction in the energy  interval length is given. The approach chosen there 
does not use the spectral shift function and thus avoids the logarithmic term.
For Anderson models on $l^2(\ZZ^d)$ with H\"older continuous 
distribution of the potential in \cite{Krishna} an alternative proof of the results in 
\cite{HundertmarkKNSV-06} was presented.

The bounds proven in this paper for quantum graph operators have a Wegner constant which is 
independent of the position of the energy interval. The same situation is encountered for 
discrete operators, whereas for alloy-type operators on $L^2(\RR^d)$ the Wegner constant 
grows with the energy. Note that this uniformity in our case is only	 possible since 
we assume a covering condition on the single site potentials. 
For site potentials of small support the uniformity does not 
hold, cf.~\cite{HelmV}.

For more background on the general theory of random Schr\"odinger operators see for instance \cite{CarmonaL-90,PasturF-92}
and for more information on Wegner estimates and the integrated density of states \cite{KirschM-07,Veselic-06b}, e.g.

The next section contains the results of the paper and the last section the proofs.

\section{Model and results}

\begin{dfn}
Let $V$ and $E$ be countable sets and $\CG$ a map
\[
\CG \colon E \to V\times V \times [0,\infty), \quad e \mapsto
(\iota(e),\tau(e),l_e).
\]
We call the triple $G=(V,E,\CG)$ a metric graph, elements of $V=V(\G)$ vertices, elements of $E=E(\G)$ edges,
$\iota(e)$ the initial vertex of $e$, $\tau(e)$ the terminal vertex of $e$ and
$l_e$ the length of $e$. Both $\iota(e)$ and $\tau(e)$ are called endvertices of $e$,
or incident to $e$. The two endvertices of an edge are allowed to coincide.
The number of edges incident to the vertex $v$ is called the
\textit{degree of $v$}. We assume that all vertices have finite degree.
\end{dfn}
\smallskip

Each edge $e$ will be identified with the open interval $\,(0,l_e)$,
where the point $0$ corresponds to the vertex $\iota(e)$ and $l_e$
to $\tau(e)$.
The identification of edges by intervals allows us to define in a
canonical way the length of a path between two points in $G$. 
Taking the infimum over the lengths of
paths connecting two given points in $G$, one obtains a distance
function $d \colon G\times G\to [0,\infty)$. Since we assumed that
each vertex of $G$ has bounded degree, the map $d$ is indeed a
metric, cf.~for instance Section 2.2 in \cite{Schubert-06}. Thus
we have turned $G$ into a metric space $(G,d)$.

For a finite subset $\Lambda \subset E$ we define the subgraph
$\G_\Lambda $ by deleting all edges $e \in E\setminus \Lambda$
and the arising isolated vertices. We denote the set of vertices of
$G_\Lambda$ by $V_\Lambda$, the set of vertices $v\in V_\Lambda$ with $\deg_{G_\Lambda}v<\deg_G v$
by $V_\Lambda^\partial$, and its complement $V_\Lambda \setminus V_\Lambda^\partial$
by $V_\Lambda^i$. Elements of $V_\Lambda^\partial$
are called boundary vertices of $G_\Lambda$
and elements of  $V_\Lambda^i$ interior vertices of $G_\Lambda$.

For any $\Lambda \subset E$ the Hilbert spaces $L_2(\G_\Lambda)$
have a natural direct sum representation
$L_2(\G_\Lambda)=\oplus_{e\in \Lambda}L_2(0,l_e)$. In particular
for $\Lambda = E$ we have $L_2(\G)= \oplus_{e\in \E}L_2(0,l_e)$,
and for $\tilde\Lambda\subset \Lambda \subset E$ we have
$L_2(\G_{\tilde\Lambda})=\oplus_{e\in
{\tilde\Lambda}}L_2(0,l_e)\subset L_2(\G_\Lambda)=\oplus_{e\in
\Lambda}L_2(0,l_e)$.

For a function $\phi\colon \G \to \CC$ and an edge $e\in E$ we denote by
$\phi_e:=\phi|_e$ its restriction to $e$ (which is identified with $(0,l_e)$).
We denote by $C(\G)$ the space of continuous, complex-valued functions on the metric space $(\G,d)$.
Similarly, $C(\G_\Lambda)$ denotes the space of continuous, complex-valued functions on the metric sub-space $(\G_\Lambda,d)$.
For each $v \in V$, any edge $e$ incident to $v$, and function $f\in W^{2,2}(e) \subset C^1(e)\cong C^1(0,l_e)$
we define the boundary value $f(v)$ by continuity and the derivatives by
\begin{align}
\partial_e f(v) &:= \partial_e f(0)   :=  \lim_{\epsilon\searrow 0} \frac{f(\epsilon)-f(0)}{\epsilon} & \text{ if } v = \iota(e)
\intertext{and}
\partial_e f(v) &:= \partial_e f(l_e) :=  \lim_{\epsilon\searrow 0} \frac{f(l_e-\epsilon)-f(l_e)}{\epsilon} & \text{ if } v = \tau(e).
\end{align}
Note that, since $f|_e\in W^{2,2}(e)$ the function $f$ is not only continuously differentiable
on the open segment $(0,l_e)$, but also its derivative has
well defined limits at both boundaries $0$ and $l_e$.
Our sign convention ensures that $\partial_e f(v)$ is the inward normal derivative,
and is independent of the orientation of the edge induced by $\iota,\tau$.

For any $\Lambda \subset E$ it will be convenient to use the
following Sobolev space
\begin{align*}
W^{2,2}(\Lambda):= \oplus_{e\in \Lambda } W^{2,2}(e) \subset
C^1(\Lambda) := \oplus_{e\in \Lambda }  C^{1}(e)
\\
\text{ with the norm }\|\phi\|^2_{W^{2,2}(\Lambda)}:= \sum_{e\in \Lambda} \|\phi_e\|^2_{W^{2,2}(0,l_e)} .
\end{align*}
Note that this space is defined on the edge set only and does not
see the graph structure of $G$. 

For  $f\in W^{2,2}(\Lambda)$ and each vertex $v$ we gather the boundary values
$f_e (v)$ over all edges $e$ adjacent to $v$ in a vector
$f(v):= \{f_e (v) : e \in E, v \text{ incident to } e\}$. 
Similarly, we gather the boundary values of $\partial_e f_e (v)$ over all
edges $e$ adjacent to $v$ in a vector $\partial f (v)$.

Given the boundary values of functions, we can now dicuss the concept
of boundary condition. Here we use material from \cite{KostrykinS-99b,Harmer-00} to which we
refer for further details and proofs.  A single-vertex boundary condition at $v\in V$
is a choice of subspace $S_v$ of $\CC^{\deg v}\times\CC^{\deg v}$ with
dimension $\deg v$ such that
 $$
 \eta((s,s'),(t,t')):=\langle s',t\rangle - \langle s, t'\rangle
 $$
vanishes for all $(s,s'), (t,t')\in S_v$.  An $f\in W^{2,2}(\Lambda)$ is said to
satisfy the single-vertex boundary condition $S_v$ at $v$ if $(f(v),\partial f(v))$ belongs
to $S_v$.
A field  of single-vertex boundary conditions $S:=\{S_v : v\in V_\Lambda\}$ will be called
boundary condition. Given such a field, we obtain a selfadjoint realization $\Delta^\Lambda$
of the Laplacian $\Delta$ on $L^2 (E^\Lambda)$ by choosing the domain
 $$
 \cD(\Delta^\Lambda) :=\{f\in W^{2,2} (\Lambda) :\forall v : (f(v),\partial f(v))\in S_v\}.
 $$

Particularly relevant boundary conditions are Dirichlet boundary
conditions with subspace $S^D$ consisting of all those $(s,s')$ with $s=0$, Neumann conditions with
subspace $S^N$ consisting of all those $(s,s')$
with $s'=0$, and Kirchhoff (also known as free or standard) boundary conditions $S^K$ consisting
of all $(s,s')$ with $s$ having all components equal and $s'$ having the sum over its components equal to $0$.

We define the linear operator
\[
-\Delta^\Lambda \colon \cD(\Delta^\Lambda) \to L^2(\Lambda)
\]
by the rule
\[
(-\Delta^\Lambda f)(x) := -\frac{\partial^2 f_e(x)}{\partial x^2 }
\]
if $x \in G^\Lambda$ is contained in the edge $e$. This way the function
$-\Delta^\Lambda f$ is defined on the set $E^\Lambda\subset G^\Lambda$, whose complement $V^\Lambda=G^\Lambda\setminus E^\Lambda$
in the metric space $G^\Lambda$ has Hausdorff measure zero.

In our application for the Wegner estimate, we will start with boundary conditions defined on the graph $G$ and then restrict
to the induced subgraph $G_\Lambda$ for a finite subset $\Lambda\subset E(G)$ as described above.
On $V_\Lambda^i$ the boundary conditions will be induced by those on $V(G)$, but on $V_\Lambda^\partial$ there is no canonical choice.
We choose to put Dirichlet conditions on $V_\Lambda^\partial$ in order to define the restriction $-\Delta^\Lambda$ of $-\Delta^{E(G)}$ unambiguously.

An \emph{alloy-type potential} is a stochastic process $\cV\colon \Omega\times G \to \RR$
of the form $\cV_\omega=\sum_{e\in E} \omega_e \, u_e$, with the conditions
outlined in the following.

 The
\emph{coupling constants} $\omega_e, e \in E$, are a sequence
of bounded random variables which are independent and identically
distributed with distribution $\mu$. We call $\mu$ the \emph{single site
distribution} since in our model each edge $e$ is the site of a single
perturbation controlled by $\omega_e$. Note that $u_e$ is not necessarily supported
on $e$, see below. The expectation of the product measure
$\PP:=\bigotimes_{e\in E} \mu$ is denoted by $\EE$. 

The family of single site potentials $u_e, e\in E$, is assumed to fulfill a covering condition and a summability condition:
\begin{dfn}
The family of single site potentials $u_e, e\in E$, is said to fulfill a \emph{covering condition}
with lower bound $\kappa>0$ if, 
for each finite set of edges $\Lambda$, there is a finite set of edges $\Lambda^u$ such that 
\[
\sum_{e\in \Lambda^u} u_e \ge \kappa 
\]
holds on the graph $G_{\Lambda}$.
\end{dfn}
For the following definition, recall that for a metric graph $\tilde G$ with finite set of edges $\tilde E$ and length function $e \mapsto l_e$
the volume is given by $\vol \tilde G=\sum_{e\in\tilde E} l_e$. In contrast to this, $|\Lambda|$ denotes the number of edges in 
$\Lambda \subset E$.
\begin{dfn}
Denote by $\Lambda_e$ the minimal set of edges containing the support of $u_e|_\Lambda$
and by $V_e^\partial$ the boundary vertices of the induced subgraph $G_{\Lambda_e}$.
(Here for simplicity we suppress the dependence of $V_e^\partial$  on $\Lambda$.)
Then, the family of single site potentials $u_e, e\in E$ is called \emph{summable} if there are constants $C_j,j=1,2,3$, such that
\begin{equation} \label{eqn:summable}
\begin{aligned}
 \sum_{e\in\Lambda^u} \sum_{v\in V_e^\partial}\deg v &\leq C_1 |\Lambda| \quad\text{(finite degree property)}, \\
 \sum_{e\in\Lambda^u} \sqrt{\|u_e\|_\infty} \vol G_{\Lambda_e} &\leq C_2 |\Lambda| \quad\text{($L^2$-boundedness)}, \\
 \sum_{e\in\Lambda^u} |\Lambda_e| &\leq C_3 |\Lambda| \quad\text{(volume growth)}
\end{aligned}
\end{equation}
for each finite set of edges $\Lambda$.
\end{dfn}
In particular, this holds if $u_e$ is supported on $e$, uniformly bounded above and away from $0$ 
(so that $\Lambda^u=\Lambda$, $\Lambda_e=e$) and there are uniform bounds on vertex degrees and edge lengths.
But our definition is much more general. For instance, decreasing edge lengths can compensate for potential growth and vice versa.

In the following we consider for a finite subset $\Lambda \subset E$ 
and an alloy-type potential whose family of single site potentials fulfills the covering condition and is summable 
the restriction $\cV_\omega^\Lambda = \cV_\omega \chi_\Lambda$.
On $\cD(\Delta^\Lambda)$ we define a  
\emph{random Schr\"odinger operator of alloy-type} by $H_\omega^\Lambda =-\Delta^\Lambda +\cV_\omega^\Lambda$.
Since the potential is bounded, $H_\omega^\Lambda $ is selfadjoint with boundary conditions described 
by the field $\{S_v : v\in V_\Lambda\}$ and lower semi-bounded. In the following we will be dealing exclusively with Schr\"odinger operators 
$H_\omega^\Lambda $ on \emph{finite} edge sets $\Lambda$.

\begin{rem}
For random Schr\"odinger operators on the whole, infinite graph $G$ one needs to impose 
more restrictive conditions if one wants to ensure that the Schr\"odinger operator is lower semi-bounded.
For instance, Dirichlet, Neumann, standard/free/Kirchhoff conditions and others with ``$L_+=0$'' lead to positive graph Laplacians (see \cite{KostrykinS-06}).
\end{rem}

Finally, we introduce the \emph{modulus of continuity} of the distribution of the single site distribution, for
$\varepsilon>0$, as
\begin{equation}
\label{definition-s-mu-epsilon}
s(\mu,\varepsilon)=\sup\{\mu([\lambda-\varepsilon,\lambda +\varepsilon]) \mid \lambda \in \RR\} .
\end{equation}

With these definitions, we can formulate
\begin{thm} \label{thm:wegner}
\label{t-WE} Let $\cV_\omega$ be an alloy-type potential. 
Then there exists a constant $C_W$ such that for all $\lambda\in\RR$, 
all finite sets of edges $\Lambda$ and all $\varepsilon \le 1/2$
\begin{equation}
\label{e-WE}
\EE\{\Tr [ \chi_{[\lambda-\varepsilon,\lambda +\varepsilon]}(H_\omega^\Lambda) ]\}
\le C_W \ s(\mu,\varepsilon) \, |\Lambda| \, .
\end{equation}
\end{thm}
\begin{rem}[integrated density of states]
For the discussion in this remark we assume that all edge lengths $l_e \equiv  l$ are equal, 
that all single site potentials $u_e$ have the same shape, that the boundary condition $S_v$ at the vertex $v\in V$
depends only on the degree $\deg (v)$, and that there is a uniform bound
$d_+ := \sup_{v \in V} \deg(v)< \infty$ on the vertex degree. 
In that case $H_\omega^\Lambda$ is 
selfadjoint and lower semi-bounded even for $\Lambda = E$, and the domain of the operator is independent of $\omega \in \Omega$.
For an energy $\lambda\in \RR$ and an exhaustion 
$\Lambda_1 \subset \Lambda_2 \subset \Lambda_3 \subset\dots$
of the edge set $E$ consider the sequence of random variables 
\[
N_\omega^n(\lambda):= \frac{1}{\vol G_{\Lambda_n}}\Tr [ \chi_{(-\infty,\lambda ]}(H_\omega^{\Lambda_n}) ].
\]

Under certain additional conditions one can show that the sequence of distribution functions $N_\omega^n(\lambda), n \in \NN$
converges and that the resulting limiting distribution function $N$ is independent of $\omega \in \Omega$ almost surely.
If this is the case Theorem \ref{thm:wegner} implies the following continuity property of $N$:
\[
\forall \varepsilon \in[0,1/2], \forall \lambda \in \RR \quad : \quad 
N(\lambda+ \epsilon) -N(\lambda) \le C_W \ s(\mu,\varepsilon)
\]
Thus the  integrated density of states inherits the continuity modulus of the single site measure $\mu$.
Let us stress that the Wegner constant $C_W$ is energy independent.  

For a particular example of a random Schr\"odinger operator on a metric graph with $\ZZ^d$-structure the construction of the 
integrated density of states has been carried out in \cite{HelmV} following the strategy of \cite{KirschM-82c}.
\end{rem}

Wegner estimates play an important role in the proof of spectral localization for random Schr\"odinger operators
via the so called multiscale analysis.
For an alloy type model on a $\ZZ^d$-metric graph localization  has been proven in \cite{ExnerHS} 
using (weaker) Wegner estimates and multiscale analysis.

\section{Proofs}
Consider a pair of selfadjoint, lower semi-bounded operators $H_1,H_2$. If the spectrum of both $H_1$ and $H_2$ is
purely discrete, the spectral shift function (SSF) $\xi(\cdot) = \xi(\cdot,H_2,H_1)$ is defined as the difference
of the eigenvalue counting functions, i.e.
\[
\xi(\lambda):= \Tr[\chi_{(-\infty,\lambda]}(H_2) -\chi_{(-\infty,\lambda]}(H_2)] .
\]
If we merely assume that the difference $H_2-H_1$ is trace class, 
then there is still a unique function $\xi$ such that \emph{Krein's trace identity}
\begin{equation}
\label{e-KTI}
\Tr \left [ \rho (H_2)-\rho (H_1)  \right  ]
= \int \rho'(\lambda) \, \xi(\lambda, H_2,H_1) \, d\lambda
\end{equation}
holds for all $\rho \in C^\infty$ with compactly supported derivative.
In the case of operators with discrete spectrum both definitions of the function $\xi$ coincide.
We can weaken the assumption on the operator pair further. Assume that there exists 
a monotone, smooth function $g\colon \RR \to [0,\infty)$ which is bounded on the spectra
of $H_1$ and $H_2$ and such that $g(H_2)-g(H_1)$  is trace class. In that case the definition
\begin{equation}
\label{e-IP}
\xi(\lambda, H_2, H_1)
:= 
\sign(g')\ \xi\big(g(\lambda), g(H_2), g(H_1)\big) .
\end{equation}
makes sense and is independent of the choice of the function $g$. 

\begin{lem}  \label{Potential-Graph}
Let ${\tilde G}$ be a finite or infinite metric graph, ${\tilde\Lambda}$ a finite subset of its edges, 
$-\Delta$ a selfadjoint realization of the Laplacian on $L^2({\tilde G})$ and 
$W_1,W_2$ two potentials acting as bounded operators on $L^2({\tilde G})$ such that $\supp (W_2-W_1) \subset \tilde G_{\tilde\Lambda}$.
Set $H_j=-\Delta +W_j,j=1,2,$ and assume that the SSF $\xi_{H_1,H_2} $ is well defined.
Denote the restriction of $H_j$ to $L^2({\tilde\Lambda})$ by $h_j, j=1,2$. Then we have
\[
\left|\xi_{H_1,H_2} (\lambda) \right| 
\le \sum_{v \in V_{\tilde\Lambda}^\partial} \deg(v) +\left|\xi_{h_1,h_2} (\lambda) \right|. 
\]
\end{lem}
\begin{proof}
The basic idea is to decouple the interior of ${\tilde G}_{\tilde\Lambda}$ from the exterior by choosing appropriate
boundary conditions on $V_{\tilde\Lambda}^\partial$; this is in the spirit of the statements around equation (9) -- (11) in \cite{HelmV}.
So, let $H_j^D$ be $H_j$ but with Dirichlet conditions in $V_{\tilde\Lambda}^\partial$.
Then $\left| \xi_{H_j, H_j^D}(\lambda)\right| \leq \sum_{v \in V_{\tilde\Lambda}^\partial} \deg(v)$
according to Corollary 11 and Lemma 13 of \cite{GruberLV}.
Now, the $H_j^D$ decompose into a direct sum of exterior and interior parts, the former coinciding by assumption,
the latter being given by $h_j$. This proves the assertion.
\end{proof}
We will apply this lemma for a set of edges ${\tilde\Lambda}:= {\Lambda}_e$ containing the support of $u_e$,
with finite $\tilde G=G_{\Lambda_e}$ such that the SSF exists automatically. 

\begin{lem}  \label{Potential-Kante}
Let $-\Delta$ be a selfadjoint realization of the Laplacian on an arbitrary  finite graph $\tilde G$  
and let  $W_1,W_2$ be bounded potentials on $L^2(\tilde G)$.
Set $H_j=-\Delta +W_j,j=1,2$.
Then 
\[
\left|\xi_{H_1,H_2} (\lambda) \right| \leq \left(\sqrt{\|W_1\|}+\sqrt{\|W_2\|}\right) \frac{\vol \tilde G}\pi  +5|E(\tilde G)|
\]
\end{lem}
where $\vol \tilde G= \sum_{e\in E(\tilde G)} l_e$ is the one-dimensional volume of $\tilde G$.
\begin{proof}
This is an extension  of Lemma 14 in \cite{GruberLV}, where this is proved for metric graphs with edges
of length $1$. Using the same proof but keeping track of the lengths $l_e$ yields the desired estimate.
\end{proof}

\begin{proof}[Proof of Theorem~\ref{thm:wegner}]
Let $\rho$ be a smooth, monotone switch function
$\rho:=\rho_{\lambda,\varepsilon}\colon \RR \to [-1,0]$. \label{p-rho} By a switch function we mean that
for a positive $\varepsilon \le 1/2$, $\rho$  has the following properties: $ \rho\equiv -1$ on
$(-\infty,\lambda-\varepsilon]$, $\rho\equiv 0$ on $[\lambda+\varepsilon,\infty)$ and $\|\rho'\|_\infty \le
1/\varepsilon$. 
Then
\[
\chi_{[\lambda-\varepsilon,\lambda +\varepsilon]} (x) \le \rho(x+2\varepsilon) -\rho(x-2\varepsilon)
\]
We may assume without loss of generality $\sum_{e\in \Lambda^u} u_e\ge 1$, i.e.\ $\kappa=1$.
By the min-max principle for eigenvalues, we conclude
\[
\Tr  [\rho (H_\omega^\Lambda+ \varepsilon)] \le \Tr\Big[\rho (H_\omega^\Lambda+ \varepsilon\sum_{e\in \Lambda^u} u_e) \Big] .
\]
Let $\Lambda^u$ be as above. Then $\sum_{e\in\Lambda^u}u_e(x)\geq1$ for $x\in G_\Lambda$.
$\Lambda^u$ contains $L:=|\Lambda^u|$ edges. We
enumerate the edges in $\Lambda^u$ by $ e\colon \{1, \dots, L\} \to \Lambda^u $, $n\mapsto e(n)$, and set
\[
W_0 \equiv 0, \quad  W_n =\sum_{m=1}^{n} u_{e(m)}, \qquad n=1,2,\dots, L
\]
Thus
\begin{align}
\chi_{[\lambda-\varepsilon,\lambda +\varepsilon]} (H_\omega^\Lambda) 
& \le       \nonumber
 \rho(H_\omega^\Lambda+2\varepsilon)-\rho(H_\omega^\Lambda-2\varepsilon)
\\      \label{e-project}
& \le
\rho(H_\omega^\Lambda-2\varepsilon+4\varepsilon W_{L})-\rho(H_\omega^\Lambda-2\varepsilon)
\\      \nonumber
& =
\sum_{n=1}^{L} \rho(H_\omega^\Lambda+2\varepsilon+4\varepsilon  W_{n})-
\rho(H_\omega^\Lambda-2\varepsilon+4\varepsilon  W_{n-1})
\end{align}
We fix $n \in \{1, \dots, L\}$, define 
\[
\omega^\perp := \{\omega_e^\perp\}_{e \in \Lambda^u}, \qquad
\omega_e^\perp :=\begin{cases} 0 \quad &\text{if }e=e(n), \\
\omega_e \quad &\text{if } e\neq e(n), \end{cases}
\]
and set
\[
\phi_{n}(\eta) = \Tr\bigl[\rho(H_{\omega^\perp}^\Lambda-2\varepsilon +4\varepsilon W_{n-1}+\eta
u_{e(n)})\bigr], \quad \eta\in\RR.
\]
The function $\phi_{n}$ is continuously differentiable, monotone increasing and bounded.
By definition of $\phi_{n}$,
\[
\Tr [ \rho(H_\omega^\Lambda-2\varepsilon+4\varepsilon  W_n))
-\rho(H_\omega^\Lambda-2\varepsilon+4\varepsilon  W_{n-1})]
=
\phi_{n}(\omega_{e(n)}+4\varepsilon)-\phi_{n}(\omega_{e(n)})]
\]
since $\phi_{n}(\eta) = \Tr\bigl[\rho(H_{\omega}^\Lambda-2\varepsilon +4\varepsilon W_{n-1}+(\eta
-\omega_{e(n)})u_{e(n)})\bigr]$,
so that
\[
\EE_{\omega_{e(n)}} \{ \Tr [ \rho(H_\omega^\Lambda-2\varepsilon+4\varepsilon  W_n))
-\rho(H_\omega^\Lambda-2\varepsilon+4\varepsilon  W_{n-1})] \}
=
\int [\phi_{n}(\omega_{e(n)}+4\varepsilon)-\phi_{n}(\omega_{e(n)})]\, d\mu(\omega_{e(n)})
\]
where $\EE_{\omega_{e(n)}}$ denotes the expectation with respect to the random variable $\omega_{e(n)}$ only.
Let $\supp(\mu)\subset(a,b)$.
Using Lemma 6 in \cite{HundertmarkKNSV-06}  
 we have
\begin{align*}
\int [\phi_{n}(\omega_{e(n)}+4\varepsilon)-\phi_{n}(\omega_{e(n)})]\, d\mu(\omega_{e(n)})
&\le
s(\mu,4\varepsilon) [\phi_{n}(b)-\phi_{n}(a)] 
\end{align*}
Denote by $\xi_{\Lambda,n}$ the SSF associated to the pair of operators $H_n(a),H_n(b)$ on $L^2(\Lambda)$ 
where $H_n(\eta)$ is given by $H_n(\eta):=H_\omega^\Lambda-2\varepsilon+4\varepsilon  W_{n-1}+ (\eta-\omega_{e(n)}) u_{e(n)} $.
Then by the Krein trace identity and the normalization of $\rho$
\begin{align*}
\phi_{n}(b)-\phi_{n}(a)= \int_a^b \rho' \ \xi_{\Lambda,n} \ d\lambda  
\le \|\xi_{\Lambda,n}\|_\infty 
\end{align*}
Let $\Lambda_e$, $u_e$, $V_{e}^\partial$ and $G_{\Lambda_e} $ be as in the definition of summable potentials.
By $\xi_{\Lambda_{e(n)},n}$ we denote the SSF associated to the pair $H_n(a),H_n(b)$,
but now considered as operators on $L^2(\Lambda_{e(n)})$.  
Apply Lemma \ref{Potential-Graph} to obtain:
\begin{align*}
\|\xi_{\Lambda,n}\|_\infty 
\le \sum_{v \in V_{e(n)}^\partial} \deg(v) +\|\xi_{\Lambda_{e(n)},n}\|_\infty 
\end{align*}

Now apply Lemma \ref{Potential-Kante} successively $L$ times to obtain
\[
\begin{aligned}
\EE \{\Tr [\chi_{[\lambda-\varepsilon,\lambda +\varepsilon]} (H_\omega^\Lambda) ]\}
& \le s(\mu,4\varepsilon) \sum_{n=1}^{L} \left( \sum_{v \in V_{e(n)}^\partial} \deg(v) 
     + \sqrt{\|u_{e(n)}\|_\infty} \frac{\vol G_{\Lambda_{e(n)}}}\pi +  5 | \Lambda_{e(n)}| \right) \\
&\le s(\mu,4\varepsilon) (C_1+C_2/\pi+5C_3) |\Lambda|
\end{aligned}
\]
by the summability condition on the family of single site potentials.
\end{proof}

\bibliographystyle{habbrv}
\bibliography{gv}	

\begin{thebibliography}{10}

\bibitem{CarmonaL-90}
R.~Carmona and J.~Lacroix.
\newblock {\em Spectral Theory of Random {Schr\"odinger} Operators}.
\newblock Birkh\"auser, Boston, 1990.

\bibitem{CombesHK}
J.-M. Combes, P.~D. Hislop, and F.~Klopp.
\newblock An optimal {Wegner} estimate and its application to the global
  continuity of the integrated density of states for random {Schr\"odinger}
  operators.
\newblock 2006, arXiv:math-ph/0605029.

\bibitem{ExnerHS}
P.~Exner, M.~Helm, and P.~Stollmann.
\newblock Localization on a quantum graph with a random potential on the edges.
\newblock 2006, arXiv:math-ph/0612087.

\bibitem{GruberLV}
M.~J. Gruber, D.~Lenz, and I.~Veseli\'c.
\newblock Uniform existence of the integrated density of states for random
  {S}chr\"odinger operators on metric graphs over {$\mathbb{Z}^d$}.
\newblock 2006, arXiv:math.SP/0612743.

\bibitem{Harmer-00}
M.~Harmer.
\newblock Hermitian symplectic geometry and extension theory.
\newblock {\em J. Phys. A}, 33(50):9193--9203, 2000.

\bibitem{HelmV}
M.~Helm and I.~Veseli{\'c}.
\newblock A linear {Wegner} estimate for alloy type {Schr\"odinger} operators
  on metric graphs.
\newblock 2006, arXiv:math.SP/0611609.

\bibitem{HundertmarkKNSV-06}
D.~Hundertmark, R.~Killip, S.~Nakamura, P.~Stollmann, and I.~Veseli{\'c}.
\newblock Bounds on the spectral shift function and the density of states.
\newblock {\em Comm. Math. Phys.}, 262(2):489--503, 2006.

\bibitem{KirschM-82c}
W.~Kirsch and F.~Martinelli.
\newblock On the density of states of {Schr\"odinger} operators with a random
  potential.
\newblock {\em J. Phys. A: Math. Gen.}, 15:2139--2156, 1982.

\bibitem{KirschM-07}
W.~Kirsch and B.~Metzger.
\newblock The integrated density of states for random {Schr\"odinger}
  operators.
\newblock In {\em Spectral Theory and Mathematical Physics: A Festschrift in
  Honor of Barry Simon's 60th Birthday}, volume~76 of {\em Proceedings of
  Symposia in Pure Mathematics}, pages 649--698. AMS, 2007,
  arXiv:math-ph/0608066.

\bibitem{KostrykinS-99b}
V.~Kostrykin and R.~Schrader.
\newblock Kirchhoff's rule for quantum wires.
\newblock {\em J. Phys. A}, 32(4):595--630, 1999.

\bibitem{KostrykinS-06}
V.~Kostrykin and R.~Schrader.
\newblock Laplacians on metric graphs: eigenvalues, resolvents and semigroups.
\newblock In {\em Quantum graphs and their applications}, volume 415 of {\em
  Contemp. Math.}, pages 201--225. Amer. Math. Soc., Providence, RI, 2006.

\bibitem{Krishna}
M.~Krishna.
\newblock Continuity of integrated density of states -- independent randomness.
\newblock 2006, arXiv:math-ph/0609040.

\bibitem{PasturF-92}
L.~A. Pastur and A.~L. Figotin.
\newblock {\em Spectra of Random and Almost-Periodic Operators}.
\newblock Springer Verlag, Berlin, 1992.

\bibitem{Schubert-06}
C.~Schubert.
\newblock Laplace-{Operatoren} auf {Quantengraphen}.
\newblock Diplomarbeit, TU Chemnitz, 2006.

\bibitem{Veselic-06b}
I.~Veseli\'c.
\newblock Existence and regularity properties of the integrated density of
  states of random {Schr\"odinger} operators, 2006.
\newblock Habilitation thesis. TU Chemnitz. To be published in the Springer LNM
  series.

\bibitem{Wegner-81}
F.~Wegner.
\newblock Bounds on the density of states in disordered systems.
\newblock {\em Z. Phys. B}, 44(1-2):9--15, 1981.

\end{thebibliography}

\end{document}